# STABLE MARKED POINT PROCESSES

### By Tucker McElroy and Dimitris N. Politis

### U.S. Bureau of the Census and University of California, San Diego

In many contexts such as queuing theory, spatial statistics, geostatistics and meteorology, data are observed at irregular spatial positions. One model of this situation involves considering the observation points as generated by a Poisson process. Under this assumption, we study the limit behavior of the partial sums of the marked point process $\{(t_i, X(t_i))\}$, where $X(t)$ is a stationary random field and the points $t_i$ are generated from an independent Poisson random measure $\mathbb{N}$ on $\mathbb{R}^d$. We define the sample mean and sample variance statistics and determine their joint asymptotic behavior in a heavy-tailed setting, thus extending some finite variance results of Karr [*Adv. in Appl. Probab.* **18** (1986) 406–422]. New results on subsampling in the context of a marked point process are also presented, with the application of forming a confidence interval for the unknown mean under an unknown degree of heavy tails.

**1. Introduction.** Random field data arise in diverse areas such as spatial statistics, geostatistics and meteorology, to name but a few. It often happens that the observation locations of the data are irregularly spaced, this being a serious deviation from the typical formulation of random field theory, where data are located at lattice points. One effective way of modeling the observation points is by means of a Poisson process. In [3, 4], the statistical problem of mean estimation is addressed given a marked point process structure of the data where the observation locations are governed by a Poisson random measure $\mathbb{N}$ assumed to be independent of the distribution of the stationary random field itself. Karr [4] obtained central limit theorem results for the sample mean in this context, under a finite second moment assumption, and showed that the limiting variance depends on the integrated autocovariance function. The paper at hand is a first analysis of infinite-variance marked point processes.

Within the literature on dependent, heavy-tailed stationary time series, the discrete-time stochastic process

$$(1) \qquad\qquad X(t) = \int_{\mathbb{R}} \psi(t+x) \mathbb{M}(dx)$$

has been studied in [9]. In this work, $t$ is the integer index of the discrete-time process, $\mathbb{M}$ is an $\alpha$-stable random measure with Lebesgue control measure and $\psi$ is a sufficiently regular real-valued function. This is an excellent model for the types of data discussed above because $X(t)$ can be defined for any $t$ in $\mathbb{R}^d$ when $\psi$ and $\mathbb{M}$ are also extended to $\mathbb{R}^d$. A marked point process can be defined using model (1); data from a marked point process are of the form $\{t_i, X(t_i)\}$ for $i = 1, 2, \ldots$, where the points $t_1, t_2, \ldots$ are generated by a Poisson random measure. We focus on the non-Gaussian case where $\alpha < 2$ so that the variance of $X$ is infinite. Nevertheless, we study the sample variance statistic since it forms a suitable studentization of the mean; see [5]. The sample variance (as well as its square root, the sample standard deviation) is always well defined, even though the true variance may be infinite.

The second section of this paper develops some theory on continuous-time stable processes and the convergence of integrated partial sums of the data is established. This is relevant to our main discussion since a continuously observed process must be defined before we can conceive of a marked point process, the random field having to be well defined at every observation point $t$. Since the limit results for continuous-time stable processes are new and helpful for our marked point process problem, they are also included.

In the third section, we describe the marked point process situation and derive the joint asymptotics for the sample mean and sample variance. As expected, the limit is stochastic, but its randomness only comes from the stable random measure $\mathbb{M}$, not from the Poisson random measure $\mathbb{N}$. It will be seen that our results generalize the limit theory of Karr [4] to the case of infinite variance.

Finally, it is well known that subsampling is applicable in the context of a marked point process under some conditions; see Chapter 6 of [8], hereafter denoted PRW [8]. Our limit theorems permit us to verify the subsampling requirements, so a valid confidence interval for the mean is constructed in Section 4. Two methods are presented, one based on the asymptotics of the sample mean when $\alpha$ is known and one based on the asymptotics of the self-normalized sample mean when $\alpha$ is not known. These methods are tested and compared by means of a simulation study in Section 5. All technical proofs are in Appendix A.

**2. Continuous parameter processes.** In this section, we develop a limit theory for the sample mean and sample variance of an $\alpha$-stable continuous-time random field $\{X(t), t \in \mathbb{R}^d\}$. Consider an $\alpha$-stable random measure $\mathbb{M}$ with skewness intensity $\beta(\cdot)$ and Lebesgue control measure (denoted by $\lambda$) defined on the space $\mathbb{R}^d$. The random measure is independently scattered and for Lebesgue-measurable sets $A$, the distribution of $\mathbb{M}(A)$ is $\alpha$-stable with scale $\lambda(A)^{1/\alpha}$, skewness $\int_A \beta(x)\lambda(dx)/\lambda(A)$ and location 0; see [11], page 118 for details. This choice of control measure reflects our desire that the process be strictly stationary; translation invariance is a necessary condition for stationarity in such models. In addition, it is necessary for $\beta$ to have period "zero" for stationarity, that is, the skewness intensity $\beta(x)$ is constant as a function of $x$. We will denote this constant by $\beta$. Let $\psi$ be a *filter function* in $\mathbb{L}_\delta := \{f : \|f\|_\delta^\delta := \int_{\mathbb{R}^d} |f(x)|^\delta \lambda(dx) < \infty\}$ that is continuous and bounded for almost every $x$ with respect to Lebesgue measure. Then we may construct the following stochastic integral with respect to an $\alpha$-stable random measure $\mathbb{M}$ (see [11]):

$$(2) \qquad X(t) = \int_{\mathbb{R}^d} \psi(x+t)\mathbb{M}(dx),$$

where $t \in \mathbb{R}^d$. We stipulate that the number $\delta$ is in $(0, \alpha] \cap [0, 1]$ (later, we require $\psi$ to be integrable); $\alpha$ will be fixed throughout the discussion. Note that $\alpha = 2$ corresponds to a Gaussian stochastic process and has been extensively studied. It is a known fact that the random variable $X(t)$ is $\alpha$-stable with scale parameter

$$\sigma_\psi = \left( \int_{\mathbb{R}^d} |\psi(x+t)|^\alpha \lambda(dx) \right)^{\frac{1}{\alpha}} = \left( \int_{\mathbb{R}^d} |\psi(x)|^\alpha \lambda(dx) \right)^{\frac{1}{\alpha}}$$

and skewness parameter

$$\beta_\psi = \frac{\int_{\mathbb{R}^d} \psi(x+t)^{\langle\alpha\rangle} \beta(x)\lambda(dx)}{\sigma_\psi^\alpha} = \beta \frac{\int_{\mathbb{R}^d} \psi(x)^{\langle\alpha\rangle} \lambda(dx)}{\sigma_\psi^\alpha},$$

where $b^{\langle\alpha\rangle} = \text{sign}(b)|b|^\alpha$. The location parameter is zero unless $\alpha = 1$, in which case it is

$$\mu_\psi = -\frac{2}{\pi} \int_{\mathbb{R}^d} \psi(x+t)\beta(x) \log|\psi(x+t)|\lambda(dx)$$

$$= -\frac{2\beta}{\pi} \int_{\mathbb{R}^d} \psi(x) \log|\psi(x)|\lambda(dx).$$

Intuitively, we may think of $X$ as the convolution of $\psi$ and $\mathbb{M}$, in analogy with the infinite order moving average of classical time series analysis. Essentially, one runs the independent-increments $\alpha$-stable measure $\mathbb{M}$ through the linear filter $\psi$ and the resulting time series is strictly stationary with nontrivial dependence; two random variables $X(t)$ and $X(t+k)$ are independent if and only if the lag $k$ exceeds the diameter of $\psi$'s support. Of course, $\psi$ need not be compactly supported, in which case all of the variables are dependent. Hence, this construction makes for an interesting and relevant linear heavy-tailed model. As shown in Proposition 1 of [6], the model defined by (1) is well defined and stationary. That is, for each $t$, the random variable $X(t)$ is $\alpha$-stable with location zero (unless $\alpha = 1$, in which case the location is $\mu_\psi$), constant skewness and constant scale.

We will be interested in the asymptotic distribution of the partial sums. In order to consider fairly general, nonrectangular regions, we let $K$ be a "prototype region," that is, a Lebesgue-measurable set in $\mathbb{R}^d$ with measurable boundary $\partial K$ such that $\lambda(\partial K) = 0$; see [7] for background on this concept. Then let $K_n = n \cdot K$, which essentially scales the prototype region by the integer $n$. Our statistics are computed over $K_n$ and our asymptotic results are achieved as $n \to \infty$. This device allows us to consider nonrectangular regions and, at the same time, is a realistic construction. The appropriate rate of convergence for the partial sums will then be $n^{\frac{d}{\alpha}}$, as shown in Theorem 1 below. We wish to average $X(t)$ over *all* points in $K_n$, so to that end we must calculate

$$\int_{K_n} X(t)\lambda(dt). \tag{3}$$

Now a discrete sum of the random field is certainly well defined by the linearity of the $\alpha$-stable random integral, but it is not a priori clear that (3) makes sense. Thus, we introduce the following definition.

DEFINITION 1. By expression (3), we mean the limit in probability as $m \to \infty$ of

$$\sum_{t_i \in K_n^m} X(t_i)\Delta t_i, \tag{4}$$

where $K_n^m$ is a mesh of $m$ points $t_i$ in $K_n$, $\Delta t_i$ is the $d\lambda$ volume of the elements of the mesh and the mesh gets progressively finer as $m$ is increased (this is the usual 'Riemann sums' construction).

Let us establish that this definition makes sense. By using the linearity of the stable integral, (4) becomes

$$\int_{\mathbb{R}^d} \sum_{t_i \in K^m} \psi(t_i + x)\Delta t_i \mathbb{M}(dx) = \int_{\mathbb{R}^d} F_m(x)\mathbb{M}(dx), \tag{5}$$

where $F_m(x) = \sum_{t_i \in K^m} \psi(t_i + x)\Delta t_i$. Now, it follows that the limit in probability as $m \to \infty$ of (5) is $\int_{\mathbb{R}^d} F(x)\mathbb{M}(dx)$ for $F(x) := \int_{K_n} \psi(t + x)\lambda(dt)$,

provided

$$\int_{\mathbb{R}^d} |F_m(x) - F(x)|^\alpha \lambda(dx) \to 0$$

as $m \to \infty$. Since the integrands are bounded in $\mathbb{L}^1$, we may apply the Lebesgue dominated convergence theorem and obtain our result since $F_m(x) \to F(x)$ pointwise.

Thus, we have established that expression (3) makes sense and also that it is equal to

$$(6) \qquad \int_{\mathbb{R}^d} F(x)\mathbb{M}(dx) = \int_{\mathbb{R}^d} \int_{K_n} \psi(t+x)\lambda(dt)\mathbb{M}(dx).$$

In a similar fashion, one can define the integral of the second moment,

$$(7) \qquad \int_{K_n} X^2(t)\lambda(dt),$$

as a limit in probability of a Riemann sum of $X^2(t)$; unfortunately, due to squaring, a nice representation such as (6) is not possible. However, the Laplace transform of (7) is closely related to the Fourier transform of

$$(8) \qquad \int_{\mathbb{R}^d} \int_{K_n} \psi(t+x)\mathbb{B}(dt)\mathbb{M}(dx),$$

where $\mathbb{B}$ is an independent Gaussian random measure. This is shown in the proof of Theorem 1. Note that, conditional on $\mathbb{B}$, (8) is an $\alpha$-stable random variable with scale

$$\left( \int_{\mathbb{R}^d} \left| \int_{K_n} \psi(x+t)\mathbb{B}(dt) \right|^\alpha \lambda(dx) \right)^{1/\alpha}.$$

Now, we are interested in the continuous versions of the sample mean and sample variance; it is sufficient to examine the asymptotics of

$$\left( n^{-\frac{d}{\alpha}} \int_{K_n} X(t)\lambda(dt), n^{-\frac{2d}{\alpha}} \int_{K_n} X^2(t)\lambda(dt) \right),$$

which we explore through the joint Fourier–Laplace transform. For two random variables $A$ and $B \geq 0$, this is defined by

$$\phi_{A,B}(\theta, \gamma) = \mathbb{E} \exp\{i\theta A - \gamma B\}, \qquad \theta \in \mathbb{R}, \gamma \geq 0.$$

Like the joint characteristic function, the pointwise convergence of $\phi_{A_n, B_n}$ to a function which is continuous at $(0, 0)$ establishes joint weak convergence of $A_n$ and $B_n$; see [2]. The next theorem gives a complete answer to our inquiry.

THEOREM 1. *Consider a random field defined by the model given by* (1), *where* $0 < \alpha \leq 2$. *Then the sample first and second moments jointly have limit*

$$(9) \quad \left( n^{-\frac{d}{\alpha}} \int_{K_n} X(t) \lambda(dt), n^{-\frac{2d}{\alpha}} \int_{K_n} X^2(t) \lambda(dt) \right) \overset{\mathcal{L}}{\Longrightarrow} (S_\infty(\alpha), U_\infty(\alpha))$$

*as* $n \to \infty$. *Here,* $U_\infty(\alpha)$ *is nondegenerate only if* $\alpha < 2$. $S_\infty(\alpha)$ *is an* $\alpha$-*stable random variable with scale parameter* $|\Phi|$, *where* $\Phi = \int_{\mathbb{R}^d} \psi(x)\lambda(dx)$, *and skewness parameter* $\beta \cdot \text{sign}(\Phi)$. *If* $\alpha \neq 1$, *the location parameter is zero; otherwise, it is* $-\frac{2\beta}{\pi}\Phi \log|\Phi|$. *When either* $\alpha \neq 1$ *or* $\alpha = 1$ *and* $\beta = 0$, *we may write* $S_\infty(\alpha) \overset{\mathcal{L}}{=} \lambda(K)^{1/\alpha}\Phi \cdot \mathbb{M}(B)$ *for the marginal distribution. For* $\alpha < 2$, $U_\infty(\alpha)$ *is an* $\alpha/2$-*stable random variable with scale parameter*

$$2[\lambda(K)\mathbb{E}|G|^\alpha]^{2/\alpha} \int_{\mathbb{R}^d} \psi^2(x)\lambda(dx)(\cos(\pi\alpha/4))^{2/\alpha},$$

*skewness* 1 *and location* 0, *where* $G$ *is a standard normal random variable. If* $\alpha = 2$, *then* $U_\infty(\alpha)$ *is a point mass at the second moment of* $X$, *which is* $2\int_{\mathbb{R}^d} \psi^2(s)\lambda(ds)$. *We may write the marginal distribution as*

$$U_\infty(\alpha) \overset{\mathcal{L}}{=} 2[\lambda(K)\mathbb{E}|G|^\alpha]^{2/\alpha} \int_{\mathbb{R}^d} \psi^2(x)\lambda(dx)(\cos(\pi\alpha/4))^{2/\alpha} \cdot \epsilon(\alpha),$$

*where* $\epsilon(\alpha)$ *is an* $\alpha/2$-*stable subordinator. The limiting joint Fourier–Laplace transform* $\mathbb{E}\exp\{i\theta S_\infty(\alpha) - \gamma U_\infty(\alpha)\}$ *for* $\alpha < 2$ *is* [*letting* $\Phi_2 = \sqrt{\int_{\mathbb{R}^d} \psi^2(x)\lambda(dx)}$]

$$(10) \quad \begin{aligned} &\exp\Bigg\{ -\lambda(K)\mathbb{E}|\theta\Phi + \sqrt{2\gamma}\Phi_2 G|^\alpha \left( 1 - i\beta \frac{\mathbb{E}(\theta\Phi + \sqrt{2\gamma}\Phi_2 G)^{\langle\alpha\rangle}}{\mathbb{E}|\theta\Phi + \sqrt{2\gamma}\Phi_2 G|^\alpha} \right) \\ &\quad - i\lambda(K)\frac{2\beta}{\pi}1_{\alpha=1}\mathbb{E}\Big[ \left(\theta\Phi + \sqrt{2\gamma}\Phi_2 G\right) \log\Big|\theta\Phi + \sqrt{2\gamma}\Phi_2 G\Big| \Big] \Bigg\} \end{aligned}$$

*for* $\theta$ *real and* $\gamma \geq 0$.

REMARK 1. One can develop confidence intervals via this result. However, due to considerations of space, we will give the details only in the marked point case, the continuous case being much simpler.

**3. Marked point processes.** We will now consider the more intricate situation wherein the observation locations of the random field are themselves random. It often happens in statistical problems that random field data is not observed at lattice points, but instead at points scattered around the observation region with no discernible pattern. Frequently, we can model this situation through the employment of a random measure for the point locations. Generally, this probabilistic structure is referred to as a *marked point process* $\tilde{N}$:

$$\tilde{N} = \sum_i \epsilon_{(T_i, X(T_i))}$$

for $\epsilon_x(A) = 1$ if $x \in A$ and equals 0 otherwise. See [4] for a treatment of marked point processes for $\mathbb{L}_2$ random fields. If we wish to impose the condition that the distribution of points does not depend on the location of the observation region, only on its size and shape, then we say that the random measure is *spatially homogeneous*—this is similar to a stationarity assumption. Also, it is often sensible to assume that the distribution of points in one observation region is independent of the distribution of points in another disjoint observation region—the "independent scattering" property. It turns out that a homogeneous Poisson random measure (PRM) satisfies these properties, and is therefore a reasonable model in many scenarios. So, let $\mathbb{N}$ denote a PRM with mean measure $\Lambda$. $\mathbb{N}$ is sometimes denoted $PRM(\Lambda)$, as explained in the following:

DEFINITION 2. We say that $\mathbb{N}$ is $PRM(\Lambda)$ on the measure space $\{\mathbb{R}^d, \mathcal{B}, \Lambda\}$ (where $\mathcal{B}$ consists of the Borel sets in $\mathbb{R}^d$) if and only if it is an independently scattered, countably additive random measure that satisfies

(11) $$\mathbb{P}[\mathbb{N}(A) = k] = \exp\{-\Lambda(A)\} \frac{\Lambda(A)^k}{k!}, \qquad k = 0, 1, \ldots$$

for every $A \in \mathcal{B}_0 := \{B \subset \mathcal{B} : \Lambda(A) < \infty\}$—in other words, $\mathbb{N}(A) \sim \mathcal{P}ois(\Lambda(A))$. We call $\Lambda$ the mean measure of $\mathbb{N}$.

REMARK 2. More generally, we could just define the mean measure to be $\Lambda(\cdot) := \mathbb{E}\mathbb{N}(\cdot)$. Now, if we impose the condition that $\Lambda$ be a translation-invariant measure on $\mathbb{R}^d$, then spatial homogeneity follows immediately from (11). Of course, $\Lambda$ must then be Lebesgue measure (denoted by $\lambda$, as in the previous section) modulo some constant positive multiplicative factor, that is, $\Lambda = r\lambda$ for some $r \in \mathbb{R}^+$.

We are interested in investigating the limit behavior of the sample mean and sample variance over the observation region $K_n$ (we preserve the notation from the previous section), where the data locations are now determined by the random measure $\mathbb{N}$, which is independent of the stochastic process (1). Thus, we wish to study the joint convergence of

(12) $$\left( \mathbb{N}(K_n)^{-\frac{1}{\alpha}} \int_{K_n} X(t) \mathbb{N}(dt), \mathbb{N}(K_n)^{-\frac{2}{\alpha}} \int_{K_n} X^2(t) \mathbb{N}(dt) \right)$$

as $n \to \infty$. Note that $\mathbb{N}(K_n)$ is the actual observed sample size.

THEOREM 2. *Consider a continuous-parameter random field generated from the model given by* (1), *where* $0 < \alpha \leq 2$ *and the skewness intensity* $\beta$ *is constant. Suppose that a PRM* $\mathbb{N}$ *with mean measure* $\Lambda = r\lambda$, *independent of the stochastic process, governs the distribution of observation locations. If the observation region is the set* $K_n$ *and* $\alpha < 2$, *then the normalized sample mean and sample variance computed over the observation region jointly converge in distribution to an* $\alpha$-*stable random variable* $\tilde{S}_\infty(\alpha)$ *and a positive* $\alpha/2$-*stable random variable* $\tilde{U}_\infty(\alpha)$ *as* $n \to \infty$:

$$(13) \quad \left( \mathbb{N}(K_n)^{-\frac{1}{\alpha}} \int_{K_n} X(t)\mathbb{N}(dt), \mathbb{N}(K_n)^{-\frac{2}{\alpha}} \int_{K_n} X^2(t)\mathbb{N}(dt) \right)$$
$$\stackrel{\mathcal{L}}{\Longrightarrow} (r^{-1/\alpha}\tilde{S}_\infty(\alpha), r^{-2/\alpha}\tilde{U}_\infty(\alpha)).$$

*The joint Fourier–Laplace transform of the limit (*$\theta$ *real and* $\gamma > 0$*) is given by*

$$\mathbb{E}\exp\{i\theta\tilde{S}_\infty(\alpha) - \gamma\tilde{U}_\infty(\alpha)\}$$
$$= \exp\left\{ -\tilde{\sigma}_\infty^\alpha(\theta, \gamma)\left( 1 - i\beta\frac{\tilde{\beta}_\infty(\theta, \gamma)}{\tilde{\sigma}_\infty^\alpha(\theta, \gamma)} \right) + i1_{\{\alpha=1\}}\tilde{\mu}_\infty(\theta, \gamma) \right\},$$

*with the parameters given by*

$$\tilde{\sigma}_\infty(\theta, \gamma) = (\mathbb{E}|g(\theta, \gamma)|^\alpha)^{1/\alpha},$$
$$\tilde{\beta}_\infty(\theta, \gamma) = \mathbb{E}(g(\theta, \gamma))^{\langle\alpha\rangle},$$
$$\tilde{\mu}_\infty(\theta, \gamma) = \frac{-2\beta}{\pi}\mathbb{E}[g(\theta, \gamma)\log g(\theta, \gamma)],$$
$$g(\theta, \gamma) = \theta\int_{\mathbb{R}^d}\psi(s)\mathbb{N}(ds) + \sqrt{2\gamma}\sqrt{\int_{\mathbb{R}^d}\psi^2(s)\mathbb{N}(ds)}G,$$

*for* $G$ *a standard normal random variable independent of* $\mathbb{N}$. *Hence,* $\tilde{S}_\infty(\alpha)$ *is* $\alpha$-*stable with scale* $\tilde{\sigma}_\infty(1, 0)$, *skewness* $\beta\frac{\tilde{\beta}_\infty(1,0)}{\tilde{\sigma}_\infty^\alpha(1,0)}$ *and location* $1_{\{\alpha=1\}}\tilde{\mu}_\infty(1, 0)$, *whereas* $\tilde{U}_\infty(\alpha)$ *is* $\alpha/2$-*stable with scale*

$$2\left( \cos(\pi\alpha/4)\mathbb{E}\left[ \int_{\mathbb{R}^d}\psi^2(s)\mathbb{N}(ds) \right]^{\alpha/2}\mathbb{E}|Z|^\alpha \right)^{2/\alpha},$$

*skewness one and location zero. If* $\alpha = 2$, *the limit of the sample variance* $\tilde{U}_\infty(\alpha)$ *is a point mass at* $2r\int_{\mathbb{R}^d}\psi^2(s)\lambda(ds)$. *The stated limit for the sample mean still holds when* $\alpha = 2$; *in this case,* $\tilde{S}_\infty(2)$ *is a Gaussian random variable.*

As a special case, let us fix $\alpha = 2$ in Theorem 2. The limiting squared scale $\tilde{\sigma}_\infty^2(\theta, 0)$ (half of the variance) of the partial sums is then

$$\mathbb{E}\left| \int_{\mathbb{R}^d}\psi(s)\mathbb{N}(ds) \right|^2 = \left( \int_{\mathbb{R}^d}\psi(s)r\lambda(ds) \right)^2 + \int_{\mathbb{R}^d}\psi^2(s)r\lambda(ds)$$

$$= r^2 \int_{\mathbb{R}^d} \int_{\mathbb{R}^d} \psi(x)\psi(x+t)\lambda(dx)\lambda(dt) + r \int_{\mathbb{R}^d} \psi^2(s)\lambda(ds)$$

$$= \frac{1}{2}r \int_{\mathbb{R}^d} \tau(t)\lambda(dt) + \frac{1}{2}r\tau(0),$$

where $\tau(t) = \mathrm{Cov}(X(t), X(0)) = 2r\int_{\mathbb{R}^d} \psi(x)\psi(x+t)\lambda(dx)$ is the covariance (or codifference) function. Therefore, the variance of the limit $r^{-1/2}\tilde{S}_\infty(2)$ is $\int_{\mathbb{R}^d} \tau(t)\lambda(dt) + \tau(0)$, which agrees with equation (3.10) of [4] in the finite-variance case.

REMARK 3. Since the codifference function $\tau(t) = 2\|\psi\|_\alpha^\alpha - \|\psi(\cdot + t) - \psi\|_\alpha^\alpha$ is a natural generalization of the covariance function for $\alpha < 2$, one may be tempted to conjecture that $\mathbb{E}|\int_{\mathbb{R}^d} \psi(s)\mathbb{N}(ds)|^\alpha = \frac{1}{2}(\int_{\mathbb{R}^d} \tau(s)\lambda(ds) + \tau(0))$. This is, in fact, false, as evaluation on a simple $\psi$ will confirm.

Corollary 1 follows from Theorem 2 by the continuous mapping theorem.

COROLLARY 1. *Under the same assumptions as Theorem 2, the self-normalized mean converges as $n \to \infty$:*

$$\frac{\int_{K_n} X(t)\mathbb{N}(dt)}{\sqrt{\int_{K_n} X^2(t)\mathbb{N}(dt)}} \overset{\mathcal{L}}{\Longrightarrow} \frac{\tilde{S}_\infty(\alpha)}{\sqrt{\tilde{U}_\infty(\alpha)}}.$$

REMARK 4. Note that no knowledge of $\alpha$ is required to compute the self-normalized mean. It is also interesting that the limit does not depend on $r$. The ratio is nonconstant, since a squared $\alpha$-stable variable never has an $\alpha/2$-stable distribution.

**4. Subsampling applications.** This section of the paper describes how subsampling methods may be used for practical application of the results in the previous section. The idea is to use the subsampling distribution estimator as an approximation of the limit distribution of our mean-centered statistic from a marked point process that depends on unknown parameters (including $\alpha$); this procedure will yield approximate quantiles for its sampling distribution and thus confidence intervals for the mean can be formed. For more details and background on these methods, see PRW [8].

In order to use subsampling, it is desirable that the model satisfy certain mixing conditions. Strong mixing is one common condition on the dependence structure which is sufficient to ensure the validity of the subsampling theorems; see [10] for its introduction. Many time series models satisfy the assumption of strong mixing; for Gaussian processes, the summability of the autocovariance function implies the strong mixing property. In the case where $d = 1$, if our process (1) is symmetric, then the strong mixing condition is always satisfied, as Proposition 3 of [6] demonstrates. The mixing condition that is needed in general is somewhat more complicated and is defined in the next subsection.

4.1. *Subsampling with known index $\alpha$.* In this subsection, assume that $\alpha$ is known and greater than 1. Consider the location-shifted model

$$(14) \qquad Z(t) := X(t) + \mu = \int_{\mathbb{R}^d} \psi(x + t) \mathbb{M}(dx) + \mu.$$

This is the appropriate model for stable stationary random fields with nonzero location. Note that since $\alpha > 1$, the location parameter of $X(t)$ is zero. For applications, we will suppose that our data are the observations $\{Z(t_i) : t_i \in K_n\}$ for some specified observation region $K_n$ and a random collection of $\mathbb{N}(K_n)$ points $t_i$ generated by the Poisson random measure $\mathbb{N}$. Our goal is to estimate the location parameter (i.e., the mean) $\mu$ with the sample mean

$$\hat{\mu}_{K_n} = \frac{1}{\mathbb{N}(K_n)} \int_{K_n} Z(t) \mathbb{N}(dt).$$

Let $\bar{\alpha}_Z(k; l_1)$ be the mixing coefficients defined in PRW ([8], page 141). Since by (13) the sample mean converges, we can apply Theorem 6.3.1 of PRW [8] to this situation. Let $K_n(1 - c) := \{y \in K_n : B + y \subset K_n\}$ for $B := cK_n$, where $c = c_n \in (0, 1)$ is a sequence that tends to zero as the diameter of $K_n$, denoted by $\delta(K_n)$, tends to infinity. We also require that $c_n \delta(K_n) \to \infty$ as $n \to \infty$. Since $\delta(K_n) = \delta(nK) = n\delta(K)$, this means that $1/c_n = o(n)$. In practice, since $K_n$ is not clearly defined by the data, one may take it to be the convex hull or rectangular hull of the observation points. Then it is simple to produce the scaled-down copy $B$ of $K_n$, once $c$ is chosen. Let $\hat{\mu}_{K_n, B, y}$ be the statistic $\hat{\mu}$ evaluated on the set $B + y$ for any $y \in K_n(1 - c)$ and let $L_{K_n, B}(x)$ be the subsampling distribution estimator (6.8) of PRW [8]. We will assume (6.9) of PRW [8] as a condition on the mixing coefficients; this condition is easily satisfied in the special case $d = 1$ if the random field is strong mixing (which is always true for symmetric one-dimensional stable integrals, by Proposition 3 of [6]). Then by Theorem 6.3.1 of PRW [8]

$$L_{K_n, B}(x) \xrightarrow{P} J(x)$$

as $n \to \infty$, where $J(x)$ is the cumulative distribution function of $r^{-1/\alpha} \tilde{S}_\infty(\alpha)$.

REMARK 5. It follows that an asymptotically correct $(1 - p)100\%$ confidence interval for $\mu$ is given by

$$[\hat{\mu} - L_{1-p/2} \mathbb{N}(K_n)^{1/\alpha - 1}, \hat{\mu} - L_{p/2} \mathbb{N}(K_n)^{1/\alpha - 1}],$$

where $L_p$ is the $p$th quantile of $L_{K_n, B}$, defined as $\inf\{x : L_{K_n, B}(x) \geq p\}$. Note that explicit knowledge of $\alpha$ is necessary for this construction.

4.2. *Subsampling with unknown index $\alpha$.* The method outlined above is often not immediately applicable because the rate of convergence $\tau_{\Lambda(K_n)}$ depends on $\alpha$, which is typically unknown; thus, in practice, it may be necessary to estimate $\alpha$. This can be done via a subsampling estimator of the rate, as discussed in [1] or PRW ([8], Chapter 8). These methods can be extended to the marked point process scenario in a straightforward fashion; details are omitted here, but may be obtained by contacting the authors. The single important difference from [1] is that the data-driven rate $\tau_{\mathbb{N}(B+y)}$ appearing in the subsampling distribution estimator must be replaced by a deterministic rate $\tau_{\Lambda(B+y)}$; one can even use Lebesgue measure instead of the mean measure $\Lambda$ since the use of logarithms in the estimator ensures that the differences in scale are irrelevant.

Alternatively, one may be able to avoid the estimation of $\alpha$ by using a self-normalized estimate of the mean, for example, one may consider dividing by the sample standard deviation as in our Corollary 1. Suppose that $\hat{\sigma}_{K_n}$ is an estimate of scale for $X(t)$. Then we may form the ratio

$$\tau_{\mathbb{N}(K_n)} \frac{(\hat{\mu}_{K_n} - \mu)}{\hat{\sigma}_{K_n}},$$

where $\tau_u$ is the appropriate rate of convergence such that the ratio has a nontrivial weak limit. The goal is to self-normalize so that $\tau_u$ will be a known rate, that is, a rate that does not depend on unknown model parameters. A leading example is to self-normalize so that an asymptotic result with $\tau_u = \sqrt{u}$ holds. This is an improvement over the convergence rate of the (unnormalized) sample mean, where $\tau_u$ depends on $\alpha$, which is unknown.

We consider, generally, the scenario in which such a self-normalized convergence holds, that is, such that $\tau_u$ is a known rate. Corollary 1 furnishes an example of such a scenario; see the discussion following Remark 6 below. Now we adjust the definition of the subsampling distribution estimator accordingly:

$$(15) \quad L_{K_n,B}(x) := \lambda(K_n(1-c))^{-1} \int_{K_n(1-c)} 1_{\{\tau_{\mathbb{N}(B+y)}(\frac{\hat{\mu}_{K_n,B,y} - \hat{\mu}_{K_n}}{\hat{\sigma}_{K_n,B,y}}) \leqslant x\}} \, dy,$$

with $\hat{\sigma}_{K_n,B,y}$ defined similarly to $\hat{\mu}_{K_n,B,y}$. On a practical note, it is possible that $\mathbb{N}(B+y) = 0$; in this case, we get a point mass in $L_{K_n,B}(x)$ at zero, but this will not affect the asymptotics.

THEOREM 3. *Assume the convergences*

$$\tau_{\mathbb{N}(K_n)} \frac{(\hat{\mu}_{K_n} - \mu)}{\hat{\sigma}_{K_n}} \overset{\mathcal{L}}{\Longrightarrow} J,$$

$$(16) \quad a_{\mathbb{N}(K_n)}(\hat{\mu}_{K_n} - \mu) \overset{\mathcal{L}}{\Longrightarrow} V,$$

$$d_{\mathbb{N}(K_n)}\hat{\sigma}_{K_n} \overset{\mathcal{L}}{\Longrightarrow} W$$

*for positive $a_n$ and $d_n$ such that $\tau_n = a_n/d_n$ and where $W$ does not have positive mass at zero. Let $\mu$ be a parameter and assume that $\tau_u$ has the form (6.6) of PRW [8]. Let $c = c_n \in (0, 1)$ be such that $c_n \to 0$, but $c_n \delta(K_n) \to \infty$. Finally, assume the mixing condition (6.9) of PRW [8]. Then the following conclusions hold:*

   (i) *$L_{K_n, B}(x) \xrightarrow{P} J(x)$ for every continuity point $x$ of $J(x)$;*
   (ii) *if $J(\cdot)$ is continuous, then $\sup_x |L_{K_n, B}(x) - J(x)| \xrightarrow{P} 0$;*
   (iii) *letting*

$$c_{K_n, B}(1 - t) = \inf\{x : L_{K_n, B}(x) \geqslant 1 - t\},$$

*if $J(x)$ is continuous at $x = \inf\{x : J(x) \geqslant 1 - t\}$, then*

$$\mathbb{P}\{\tau_{\Lambda(K_n)}(\hat{\theta}_{K_n} - \theta) \leqslant c_{K_n, B}(1 - t)\} \to 1 - t.$$

*Thus, the asymptotic coverage probability of the interval $[\hat{\mu}_{K_n} - \tau_{\mathbb{N}(K_n)}^{-1} \times c_{K_n, B}(1 - p), \infty)$ is the nominal level $1 - p$.*

REMARK 6.   Since the subsampling distribution estimator involves Riemann integration, some numerical approximation must be made in its calculation; see Section 6.4 of PRW [8] for further details.

As a specific application of Theorem 3, consider the model given by (14) and define the sample standard deviation by

$$\hat{\sigma}_{K_n} = \sqrt{\frac{1}{\mathbb{N}(K_n)} \int_{K_n} (Z(t) - \hat{\mu}_{K_n})^2 \mathbb{N}(dt)}.$$

Hence, we have the convergence

$$\tag{17} \mathbb{N}(K_n)^{1/2} \frac{(\hat{\mu}_{K_n} - \mu)}{\hat{\sigma}_{K_n}} \xrightarrow{\mathcal{L}} \frac{\tilde{S}_\infty(\alpha)}{\sqrt{\tilde{U}_\infty(\alpha)}}$$

as $n \to \infty$, which follows immediately from Corollary 1. Note that when $\alpha = 2$, this convergence is valid, but $\tilde{U}_\infty(\alpha)$ is degenerate [it is equal to the variance of $X$, which is $2r \int_{\mathbb{R}^d} \psi^2(s)\lambda(ds)$]. The limiting ratio, however, is always absolutely continuous. Hence, we may apply Theorem 3 with $\tau_u = \sqrt{u}$ and obtain an asymptotic $1 - p$ confidence interval for $\mu$,

$$\left[\hat{\mu}_{K_n} - L_{1-p/2}\frac{\hat{\sigma}_{K_n}}{\sqrt{\mathbb{N}(K_n)}}, \hat{\mu}_{K_n} - L_{p/2}\frac{\hat{\sigma}_{K_n}}{\sqrt{\mathbb{N}(K_n)}}\right],$$

where $p = L_{K_n, B}(L_p)$. Note that no explicit knowledge of $\alpha$ is necessary for this construction; neither is it necessary to know the Poisson intensity $r$.

**5. Simulation studies.** We now focus on the above mean estimation scenario, with the asymptotic result (17). So, in this case $\tau_u = \sqrt{u}$ is a known rate. In this section we demonstrate the methods of Section 4 by means of several simulation studies. First, we illustrate how a stable marked point process can be simulated and then we discuss the practical implementation of the subsampling distribution estimators. Finally, our simulation results present the empirical coverage of the confidence intervals constructed via the subsampling methodology.

5.1. *Implementation.* Following Samorodnitsky and Taqqu ([11], page 149) (also see Resnick, Samorodnitsky and Xue [9]), the series representation of (1) is

$$X(t) = C_\alpha^{1/\alpha} \sum_{i \geq 1} \epsilon_i \Gamma_i^{-1/\alpha} \psi(U_i + t) q(U_i)^{-1/\alpha},$$

provided that $X(t)$ is symmetric $\alpha$-stable. In this representation,

- $\{\epsilon_i\}$ are i.i.d. Rademacher random variables,
- $\{\Gamma_i\}$ are the arrival times of a unit Poisson process (so they are sums of $i$ i.i.d. unit exponentials),
- $\{U_i\}$ are i.i.d. random variables with probability density function $q$,
- $C_\alpha$ is a positive constant defined in (1.2.9) of [11],

and all three of these sequences are independent of each other. We have freedom to select $q$, as long as it is a probability density function with support on the whole real line. In our simulations, we take $q$ to be the Cauchy density in $\mathbb{R}^d$; in practice, a heavy-tailed $q$ has proven more effective in producing realistic simulations. For simulation, we adopt the following procedure. First, fix $\alpha$ and determine the observation region $K_n = nK$, which can have a variety of shapes in $\mathbb{R}^d$. Also, let $\Lambda$ be Lebesgue measure for simplicity. If we want a sample of size $n$, we might choose $K_n$ such that $\Lambda(K_n) = n$, although this is no guarantee that we will obtain $n$ data points.

1. Simulate $\mathbb{N}(K_n)$ which is Poisson with mean rate $\Lambda(K_n)$.
2. Simulate $T_1, \ldots, T_{\mathbb{N}(K_n)}$ i.i.d. from a uniform distribution on $K_n$.
3. Simulate $\{\epsilon_i\}$, $\{\Gamma_i\}$ and $\{U_i\}$ for $i \leq I$, where $I$ is a predetermined threshold.
4. Determine a vector $v$, which is the "center" of the region $K_n$. Compute

$$X(T_j) = C_\alpha^{1/\alpha} \sum_{i=1}^{I} \epsilon_i \Gamma_i \psi(U_i + T_j - v) q(U_i)^{-1/\alpha}$$

for $j = 1, \ldots, \mathbb{N}(K_n)$.

Use of the centering constant $v$ is optional, but we found that it improves the quality of the simulation; theoretically, it merely introduces a deterministic lag and thus does not affect the distribution. One choice of $v$ is to let each component be defined by the various centroids. As mentioned in [11], simulation by this method is unwieldy because convergence in $I$ is slow. However, there is no alternative method for correlated stable random fields. By trial and error, we found that $I = 100$ gave a decent trade-off between simulation quality and speed; increasing $I$ to 1000 gave little visible improvement to the simulation, while greatly retarding the speed.

Next, we need to compute the subsampling distribution estimators given by (6.8) of PRW [8] and (15). The easiest method is to approximate the integrals via Monte Carlo. This is achieved by drawing a large number of random variables (we used 10,000) uniformly distributed on $K_n(1-c)$. One detail is that for a given simulated $y$, the number $\mathbb{N}(B+y)$ could be zero; this will create division by zero, problems for $\hat{\theta}_{K,B,y}$, so by convention the latter is set equal to zero in this case. In practice, this creates a point mass at zero in the subsampling distribution, but the effect is lessened by taking larger $c$ values. Another practical problem occurs when $\mathbb{N}(B+y) = 1$, which creates a division by zero issue for the self-normalized subsampling distribution estimator. In this case, we should have

$$\tau_{\mathbb{N}(B+y)}\left(\frac{\hat{\mu}_{K_n,B,y} - \hat{\mu}_{K_n}}{\hat{\sigma}_{K_n,B,y}}\right) = 1 \cdot (X(t^*) - \hat{\mu}_{K_n})/0 = \pm\infty,$$

where the sign depends on whether or not the data value at the single point $t^*$ exceeds the mean. In our computer code, we replace $\hat{\sigma}_{K_n,B,y}$ by a very small value in this case, so the resulting ratio will be a large positive or negative value, corresponding to whether $X(t^*)$ is greater or less than the mean.

5.2. *Results.* Our simulation study focuses on dimension $d = 2$, with $K_{100}$ given by both a square region (10 by 10) and a rectangular region (5 by 20). Using Lebesgue mean measure, this gives us samples of average size 100. The centering vector $v$ is simply given by the midpoints $(5,5)$ and $(2.5,10)$ for the square and rectangular regions, respectively. We simulated stable processes for $\alpha$ values $1.1, 1.2, \ldots, 1.9$ and a "Gaussian" filter function $\psi(x_1, x_2) = \exp\{-(x_1^2 + x_2^2)/2\}$ and investigated the block ratios $c = 0.1, 0.2, 0.3, 0.4$.

Method 1 assumes that $\alpha$ is known (see Section 4.1), whereas Method 2 uses a self-normalized subsampling distribution, as in Section 4.2. Each simulation was performed 1000 times. Tables 1–6 record the proportion of simulations for which the constructed confidence interval contained the true mean of zero (the standard errors are approximately 0.0095, 0.0069 and 0.0031 for $\alpha = 0.1$, 0.05 and 0.01 respectively). Clearly, both methods are sensitive to the choice of $c$. For Method 2, the coverage is a decreasing function of $c$ and an increasing function of $\alpha$, whereas for Method 1, the coverage is not monotonic in $c$ and is decreasing in $\alpha$. For most cases, an optimal value of $c$ for Method 2 would seem to lie between 0.2 and 0.3, based on the observed pattern that $c = 0.2$ resulted in overcoverage and $c = 0.3$ in undercoverage. In contrast, it seems that for high values of $\alpha$, no value of $c$ could be found to provide good coverage for Method 1. Neither method was particularly sensitive to the shape of the sampling region since results for the square and the rectangular region were similar. Although Method 1 had superior coverage for low $\alpha$, the overall performance of Method 2 was superior. In general, the choice of $c$ will depend on how important undercoverage and overcoverage are for a particular problem; $c$ is also sensitive to the shape of $K_n$. Finally, the coverage did improve with larger sample size, but for reasons of brevity, those results are not displayed here.

TABLE 1
*Coverage at nominal level 0.90, square region*

| $c$ | Method 1 | | | | Method 2 | | | |
|---|---|---|---|---|---|---|---|---|
| | **0.1** | **0.2** | **0.3** | **0.4** | **0.1** | **0.2** | **0.3** | **0.4** |
| $\alpha = 1.1$ | 0.926 | 0.961 | 0.954 | 0.869 | 0.983 | 0.909 | 0.699 | 0.577 |
| $\alpha = 1.2$ | 0.904 | 0.954 | 0.918 | 0.857 | 0.981 | 0.933 | 0.719 | 0.591 |
| $\alpha = 1.3$ | 0.807 | 0.928 | 0.909 | 0.810 | 0.991 | 0.942 | 0.752 | 0.638 |
| $\alpha = 1.4$ | 0.726 | 0.911 | 0.880 | 0.821 | 0.987 | 0.961 | 0.770 | 0.637 |
| $\alpha = 1.5$ | 0.622 | 0.854 | 0.863 | 0.806 | 0.995 | 0.954 | 0.780 | 0.687 |
| $\alpha = 1.6$ | 0.547 | 0.777 | 0.834 | 0.777 | 0.998 | 0.965 | 0.783 | 0.721 |
| $\alpha = 1.7$ | 0.440 | 0.761 | 0.802 | 0.740 | 0.997 | 0.974 | 0.798 | 0.716 |
| $\alpha = 1.8$ | 0.414 | 0.712 | 0.746 | 0.723 | 0.997 | 0.984 | 0.810 | 0.706 |
| $\alpha = 1.9$ | 0.366 | 0.670 | 0.739 | 0.691 | 0.999 | 0.986 | 0.839 | 0.727 |

Method 2's superior performance in simulation is interesting, since it also uses less information than Method 1 (it does not assume that $\alpha$ is known). This seems to corroborate the assertion that a data-driven normalization (such as the standard deviation) is superior in finite-sample to one based purely on a rate of convergence. If an extreme does not occur in the observed data, then normalization via a rate will overcompensate; conversely, if an unusual number of extremes (or an unusually large extreme) occurs, then the rate undercompensates. Use of the standard deviation instead will automatically adjust in an appropriate fashion, since it will be smaller in the first scenario and larger in the second.

TABLE 2
*Coverage at nominal level 0.95, square region*

| $c$ | Method 1 | | | | Method 2 | | | |
|---|---|---|---|---|---|---|---|---|
| | **0.1** | **0.2** | **0.3** | **0.4** | **0.1** | **0.2** | **0.3** | **0.4** |
| $\alpha = 1.1$ | 0.987 | 0.989 | 0.987 | 0.934 | 0.999 | 0.984 | 0.810 | 0.676 |
| $\alpha = 1.2$ | 0.987 | 0.992 | 0.959 | 0.881 | 0.996 | 0.984 | 0.829 | 0.663 |
| $\alpha = 1.3$ | 0.943 | 0.983 | 0.999 | 0.992 | 1.000 | 0.991 | 0.845 | 0.721 |
| $\alpha = 1.4$ | 0.925 | 0.964 | 0.938 | 0.886 | 0.998 | 0.990 | 0.857 | 0.701 |
| $\alpha = 1.5$ | 0.844 | 0.941 | 0.929 | 0.872 | 0.998 | 0.991 | 0.856 | 0.749 |
| $\alpha = 1.6$ | 0.763 | 0.871 | 0.894 | 0.850 | 1.000 | 0.993 | 0.873 | 0.785 |
| $\alpha = 1.7$ | 0.661 | 0.860 | 0.880 | 0.811 | 1.000 | 0.995 | 0.874 | 0.785 |
| $\alpha = 1.8$ | 0.616 | 0.808 | 0.820 | 0.794 | 1.000 | 0.998 | 0.890 | 0.762 |
| $\alpha = 1.9$ | 0.519 | 0.776 | 0.801 | 0.759 | 1.000 | 0.997 | 0.901 | 0.784 |

TABLE 3
*Coverage at nominal level 0.99, square region*

| $c$ | Method 1 | | | | Method 2 | | | |
|---|---|---|---|---|---|---|---|---|
| | **0.1** | **0.2** | **0.3** | **0.4** | **0.1** | **0.2** | **0.3** | **0.4** |
| $\alpha = 1.1$ | 0.999 | 0.999 | 0.999 | 0.981 | 1.000 | 1.000 | 0.926 | 0.784 |
| $\alpha = 1.2$ | 1.000 | 0.999 | 0.987 | 0.945 | 1.000 | 0.999 | 0.941 | 0.763 |
| $\alpha = 1.3$ | 0.995 | 0.995 | 0.999 | 0.997 | 1.000 | 1.000 | 0.946 | 0.806 |
| $\alpha = 1.4$ | 0.987 | 0.993 | 0.985 | 0.945 | 1.000 | 1.000 | 0.948 | 0.802 |
| $\alpha = 1.5$ | 0.964 | 0.978 | 0.966 | 0.925 | 1.000 | 1.000 | 0.952 | 0.831 |
| $\alpha = 1.6$ | 0.930 | 0.943 | 0.949 | 0.911 | 1.000 | 1.000 | 0.955 | 0.851 |
| $\alpha = 1.7$ | 0.862 | 0.927 | 0.922 | 0.876 | 1.000 | 1.000 | 0.949 | 0.863 |
| $\alpha = 1.8$ | 0.838 | 0.895 | 0.885 | 0.863 | 1.000 | 1.000 | 0.957 | 0.834 |
| $\alpha = 1.9$ | 0.734 | 0.873 | 0.873 | 0.827 | 1.000 | 1.000 | 0.976 | 0.845 |

TABLE 4
*Coverage at nominal level 0.90, rectangular region*

| $c$ | Method 1 | | | | Method 2 | | | |
|---|---|---|---|---|---|---|---|---|
| | **0.1** | **0.2** | **0.3** | **0.4** | **0.1** | **0.2** | **0.3** | **0.4** |
| $\alpha = 1.1$ | 0.948 | 0.965 | 0.926 | 0.874 | 0.983 | 0.861 | 0.582 | 0.492 |
| $\alpha = 1.2$ | 0.893 | 0.945 | 0.896 | 0.853 | 0.987 | 0.848 | 0.582 | 0.507 |
| $\alpha = 1.3$ | 0.835 | 0.925 | 0.883 | 0.838 | 0.992 | 0.894 | 0.637 | 0.548 |
| $\alpha = 1.4$ | 0.737 | 0.879 | 0.848 | 0.784 | 0.992 | 0.887 | 0.632 | 0.558 |
| $\alpha = 1.5$ | 0.631 | 0.840 | 0.804 | 0.764 | 0.997 | 0.921 | 0.675 | 0.560 |
| $\alpha = 1.6$ | 0.571 | 0.784 | 0.788 | 0.759 | 0.996 | 0.922 | 0.687 | 0.602 |
| $\alpha = 1.7$ | 0.472 | 0.727 | 0.731 | 0.693 | 0.998 | 0.940 | 0.707 | 0.601 |
| $\alpha = 1.8$ | 0.418 | 0.681 | 0.712 | 0.695 | 0.998 | 0.951 | 0.728 | 0.626 |
| $\alpha = 1.9$ | 0.396 | 0.645 | 0.709 | 0.672 | 0.997 | 0.954 | 0.737 | 0.651 |

TABLE 5
*Coverage at nominal level 0.95, rectangular region*

| $c$ | Method 1 | | | | Method 2 | | | |
|---|---|---|---|---|---|---|---|---|
| | **0.1** | **0.2** | **0.3** | **0.4** | **0.1** | **0.2** | **0.3** | **0.4** |
| $\alpha = 1.1$ | 0.991 | 0.990 | 0.979 | 0.932 | 0.998 | 0.960 | 0.673 | 0.558 |
| $\alpha = 1.2$ | 0.979 | 0.978 | 0.966 | 0.911 | 0.998 | 0.948 | 0.681 | 0.578 |
| $\alpha = 1.3$ | 0.966 | 0.973 | 0.943 | 0.901 | 0.999 | 0.953 | 0.716 | 0.599 |
| $\alpha = 1.4$ | 0.920 | 0.946 | 0.911 | 0.864 | 0.998 | 0.967 | 0.724 | 0.625 |
| $\alpha = 1.5$ | 0.859 | 0.909 | 0.882 | 0.833 | 1.000 | 0.960 | 0.752 | 0.625 |
| $\alpha = 1.6$ | 0.797 | 0.884 | 0.857 | 0.827 | 0.999 | 0.976 | 0.771 | 0.658 |
| $\alpha = 1.7$ | 0.719 | 0.821 | 0.812 | 0.782 | 0.999 | 0.988 | 0.795 | 0.666 |
| $\alpha = 1.8$ | 0.618 | 0.799 | 0.787 | 0.764 | 1.000 | 0.988 | 0.806 | 0.699 |
| $\alpha = 1.9$ | 0.567 | 0.751 | 0.773 | 0.745 | 1.000 | 0.990 | 0.806 | 0.696 |

## APPENDIX A

PROOF OF THEOREM 1. We will principally treat the $\alpha < 2$ case, since when $\alpha = 2$ the results are already known; see [4]. We will consider the joint Fourier–Laplace transform of

$$\left( n^{-d/\alpha} \int_{K_n} X(t)\lambda(dt), n^{-2d/\alpha} \int_{K_n} X^2(t)\lambda(dt) \right).$$

We first consider the case where the filter function has compact support in the set $L = \{x \in \mathbb{R}^d : |x_i| \leq l \,\forall i\}$. Then we can write

$$\mathbb{E} \exp\left\{ i\theta n^{-d/\alpha} \int_{K_n} X(t)\,dt - \gamma n^{-2d/\alpha} \int_{K_n} X^2(t)\,dt \right\}$$

TABLE 6
*Coverage at nominal level 0.99, rectangular region*

| $c$ | Method 1 | | | | Method 2 | | | |
|---|---|---|---|---|---|---|---|---|
| | **0.1** | **0.2** | **0.3** | **0.4** | **0.1** | **0.2** | **0.3** | **0.4** |
| $\alpha = 1.1$ | 0.999 | 1.000 | 0.996 | 0.975 | 1.000 | 0.994 | 0.818 | 0.655 |
| $\alpha = 1.2$ | 0.998 | 1.000 | 0.994 | 0.958 | 1.000 | 0.998 | 0.820 | 0.656 |
| $\alpha = 1.3$ | 1.000 | 0.996 | 0.978 | 0.952 | 1.000 | 0.996 | 0.843 | 0.685 |
| $\alpha = 1.4$ | 0.992 | 0.982 | 0.968 | 0.931 | 1.000 | 0.999 | 0.852 | 0.694 |
| $\alpha = 1.5$ | 0.963 | 0.970 | 0.942 | 0.889 | 1.000 | 0.997 | 0.861 | 0.707 |
| $\alpha = 1.6$ | 0.942 | 0.948 | 0.924 | 0.886 | 1.000 | 0.996 | 0.885 | 0.747 |
| $\alpha = 1.7$ | 0.889 | 0.912 | 0.882 | 0.852 | 1.000 | 1.000 | 0.902 | 0.736 |
| $\alpha = 1.8$ | 0.832 | 0.878 | 0.870 | 0.826 | 1.000 | 1.000 | 0.894 | 0.777 |
| $\alpha = 1.9$ | 0.794 | 0.856 | 0.860 | 0.813 | 1.000 | 1.000 | 0.899 | 0.776 |

$$= \mathbb{E} \exp \left\{ i\theta n^{-d/\alpha} \int_{K_n} X(t) \, dt + i\sqrt{2\gamma} n^{-d/\alpha} \int_{K_n} X(t) \mathbb{B}(dt) \right\}$$

$$= \mathbb{E} \exp \left\{ in^{-d/\alpha} \int_{\mathbb{R}^d} \left( \int_{K_n} \psi(x+t) \mathbb{W}(dt) \right) \mathbb{M}(dx) \right\},$$

where $\mathbb{W}_t$ is a Brownian motion with drift $\theta$ and volatility $\sqrt{2\gamma}$ and $\mathbb{B}$ is a Gaussian random measure independent of $\mathbb{M}$. Conditional on this Brownian motion, the scale parameter is

$$\sigma_n = \left( \frac{1}{n^d} \int_{\mathbb{R}^d} \left| \int_{K_n} \psi(x+t) \mathbb{W}(dt) \right|^\alpha dx \right)^{1/\alpha},$$

the skewness is $\beta \beta_n / \sigma_n^\alpha$, with

$$\beta_n = \frac{1}{n^d} \int_{\mathbb{R}^d} \left( \int_{K_n} \psi(x+t) \mathbb{W}(dt) \right)^{\langle \alpha \rangle} dx,$$

and the location $1_{\{\alpha=1\}} \mu_n$, with

$$\mu_n = \frac{-2\beta}{\pi} \frac{1}{n^d} \int_{\mathbb{R}^d} \left( \int_{K_n} \psi(x+t) \mathbb{W}(dt) \right) \log \left| \int_{K_n} \psi(x+t) \mathbb{W}(dt) \right| dx.$$

We only need to determine the convergence in probability of each parameter. We focus on the scale, as the other two are proved in a similar fashion. Let us define $H_x = \int_{K_n} \psi(x+t) \mathbb{W}(dt)$, which is Gaussian with mean $\theta \int_{K_n} \psi(x+t) \, dt$ and covariance

$$\mathbb{E}[(H_x - \mathbb{E}[H_x])(H_y - \mathbb{E}[H_y])] = 2\gamma \int_{K_n} \psi(x+t) \psi(y+t) \, dt.$$

Hence, all moments are bounded in $K_n$ since the variance is uniformly bounded by $2\gamma \int_{\mathbb{R}^d} \psi^2(t) \, dt$. Let $B = (0,1]^d$ and $J_j = \int_B |H_{x+j}|^\alpha \, dx$, so that

$$\sigma_n^\alpha = \frac{1}{n^d} \int_{\mathbb{R}^d} |H_x|^\alpha \, dx = \frac{1}{n^d} \sum_{j \in \mathbb{Z}^d} \int_B |H_{x+j}|^\alpha \, dx = \frac{1}{n^d} \sum_{j \in \mathbb{Z}^d} J_j,$$

where $J_j = 0$ if $x + j + t \notin L$ for all $x \in B$ and $t \in K_n$. Since the summands are bounded in probability (since, for example, the $2\alpha$ moment exists), we claim that the above expression is $o_P(1) + \frac{1}{n^d} \sum_{k \in K_n} J_{-k}$. If $J_j \neq 0$, then there must exist some $x \in B$ and $t \in K_n$ such that $j + x + t \in L$. For any set $A$, define $A(l) = \{y : d_\infty(y, A) \leq l\}$, where $d_\infty$ is the sup-norm metric on $\mathbb{R}^d$. Then $\sum_{j \in \mathbb{Z}^d} J_j = \sum_{j \in -K_n(l+1)} J_j$ since

$$j + x + t \in L \Rightarrow \max_{i=1,\dots,d} |j_i + x_i + t_i| \leq l$$

$$\Rightarrow \max_{i=1,\dots,d} |j_i + t_i| \leq l + 1 \qquad \text{for some } t \in K_n$$

$$\Rightarrow \inf_{t \in K_n} \max_{i=1,\dots,d} |j_i + t_i| \leq l + 1$$

and hence, for such a $j$,

$$d_\infty(j, -K_n) = \inf_{k \in K_n} d_\infty(j, -k) = \inf_{k \in K_n} \max_{i=1,\dots,d} |j_i + k_i| \leq l+1,$$

which implies that $j \in -K_n(l+1)$. Now, we partition into a disjoint union $-K_n(l+1) = -K_n \cup (-K_n(l+1) \setminus -K_n)$. It is simple to show that $-K_n(l+1) = -nK(l+1/n)$, where $K = K_1$, the prototype region. Hence, $-K_n(l+1) \setminus -K_n = n(-K(l+1/n) \setminus -K)$ and the count of integer lattice points in this set will be asymptotic (as $n \to \infty$) to $n^d \lambda(-K(l+1/n) \setminus -K)$, by the definition of Lebesgue measure. However, the sets in the sequence $-K(l+1/n) \setminus -K$, for any fixed $l$, are inscribed in one another. So by continuity from above,

$$\lambda(-K(l+1/n) \setminus -K) \to \lambda\left(\bigcap_{n \geq 1} -K(l+1/n) \setminus -K\right) \leq \lambda(\partial(-K)) = 0$$

since the intersection is a subset of the boundary. This shows that asymptotically, the only terms that count in $\sum_j J_j$ are those in $-K_n$. Now, these $J_{-k}$ variables are weakly dependent since the Gaussian variables $H_x$ and $H_y$ are independent if $\min_i |x_i - y_i| > l$. Hence, the law of large numbers holds and

$$\sigma_n^\alpha = o_P(1) + \frac{1}{n^d} \sum_{k \in K_n} \mathbb{E}[J_{-k}].$$

It is necessary to determine $\mathbb{E}[J_{-k}]$; here, we follow the argument given in the proof of Theorem 2 below, which is actually more complicated. Defining $\tilde{H}_x = \int_{\mathbb{R}^d} \psi(x+t)\mathbb{W}(dt)$, one can show that the difference between $\mathbb{E}|H_{x-k}|^\alpha$ and $\mathbb{E}|\tilde{H}_{x-k}|^\alpha$ tends to zero. Now,

$$\mathbb{E}|\tilde{H}_{x-k}|^\alpha = \mathbb{E}\left|\theta \int_{\mathbb{R}^d} \psi(t)\,dt + \sqrt{2\gamma}\sqrt{\int_{\mathbb{R}^d} \psi^2(t)\,dt}\, G\right|^\alpha,$$

which no longer depends on $x$ or $k$ (and $G$ denotes a standard normal random variable). With similar results for the other parameters, the joint Fourier–Laplace transform satisfies

$$\mathbb{E}\exp\left\{-\sigma_N^\alpha\left(1 - i\beta\frac{\beta_N}{\sigma_N^\alpha}\right) + i\mu_N \mathbb{1}_{\{\alpha=1\}}\right\}$$

$$\to \mathbb{E}\exp\left\{-\lambda(K)\sigma_\infty^\alpha(\theta,\gamma)\left(1 - i\beta\frac{\beta_\infty(\theta,\gamma)}{\sigma_\infty^\alpha(\theta,\gamma)}\right) + i\lambda(K)\mu_\infty(\theta,\gamma)\mathbb{1}_{\{\alpha=1\}}\right\},$$

where

$$\sigma_\infty(\theta,\gamma) = \left(\mathbb{E}\left|\theta\Phi + \sqrt{2\gamma}\Phi_2 G\right|^\alpha\right)^{1/\alpha},$$

$$\beta_\infty(\theta, \gamma) = \mathbb{E}\left(\theta\Phi + \sqrt{2\gamma}\Phi_2 G\right)^{\langle\alpha\rangle},$$

$$\mu_\infty(\theta, \gamma) = -\frac{2\beta}{\pi}\mathbb{E}\left[\left(\theta\Phi + \sqrt{2\gamma}\Phi_2 G\right)\log\left(\theta\Phi + \sqrt{2\gamma}\Phi_2 G\right)\right],$$

with $\Phi = \int_{\mathbb{R}^d}\psi(s)\lambda(ds)$ and $\Phi_2 = \sqrt{\int_{\mathbb{R}^d}\psi^2(s)\lambda(ds)}$. This convergence follows from the convergence in probability of the parameters because the exponential function in the transform is bounded. The existence of $S_\infty(\alpha)$ and $U_\infty(\alpha)$ follows from the continuity of the limiting Fourier–Laplace transform at $(0,0)$. Letting $\gamma = 0$, one recognizes the Fourier transform of an $\alpha$-stable random variable with parameters as described in Theorem 1, which is $S_\infty(\alpha)$; if $\theta = 0$, then one obtains the Laplace transform of a positive $\alpha/2$-stable random variable $U_\infty(\alpha)$, whose scale parameter is calculated in the statement of the theorem.

We remark that it is sufficient to consider convergence of the joint Fourier–Laplace transform as shown by Fitzsimmons and McElroy [2]. Finally, we must remove the truncation of the model. This is similar to (and actually easier than) the argument presented in the proof of Theorem 2 and is not repeated here. Thus, the proof is complete. □

PROOF OF THEOREM 2. Note that since $\mathbb{N}$ is a PRM, it follows from the law of large numbers, independent scattering, spatial homogeneity and shift invariance of $\Lambda$ that $\frac{\mathbb{N}(K_n)}{\Lambda(K_n)} \xrightarrow{\text{a.s.}} 1$ as $n$ tends to infinity. Thus, it suffices to examine the limit behavior of

$$(18) \qquad \left(n^{-\frac{d}{\alpha}}\int_{K_n}X(t)\mathbb{N}(dt), n^{-\frac{2d}{\alpha}}\int_{K_n}X^2(t)\mathbb{N}(dt)\right)$$

since $\Lambda(K_n) = rn^d\lambda(K)$; in the end, we must multiply our results by $r^{-1/\alpha}\times\lambda(K)^{-1/\alpha}$. Let us first consider a filter function $\psi$ with compact support in the set $L = \{x \in \mathbb{R}^d : |x_i| \leq l \,\forall i\}$. Then we can write

$$\mathbb{E}\exp\left\{i\theta n^{-d/\alpha}\int_{K_n}X(t)\mathbb{N}(dt) - \gamma n^{-2d/\alpha}\int_{K_n}X^2(t)\mathbb{N}(dt)\right\}$$

$$= \mathbb{E}\exp\left\{i\theta n^{-d/\alpha}\int_{K_n}X(t)\mathbb{N}(dt) + i\sqrt{2\gamma}n^{-d/\alpha}\int_{K_n}X(t)G(t)\mathbb{N}(dt)\right\}$$

$$= \mathbb{E}\exp\left\{in^{-d/\alpha}\int_{\mathbb{R}^d}\left(\theta\int_{K_n}\psi(x+t)\mathbb{N}(dt)\right.\right.$$

$$\left.\left. + \sqrt{2\gamma}\int_{K_n}\psi(x+t)G(t)\mathbb{N}(dt)\right)\mathbb{M}(dx)\right\}$$

by introducing a process of i.i.d. standard normal random variables $\{G(t)\}$ that are independent of $\mathbb{M}$. Conditional on $\mathbb{N}$ and the $G(t)$'s, this is an

$\alpha$-stable random variable with scale

$$\sigma_{\mathbb{N}} = \left(\frac{1}{n^d} \int_{\mathbb{R}^d} \left|\theta \int_{K_n} \psi(x+t)\mathbb{N}(dt) + \sqrt{2\gamma} \int_{K_n} \psi(x+t)G(t)\mathbb{N}(dt)\right|^{\alpha} dx\right)^{1/\alpha},$$

skewness $\beta\beta_{\mathbb{N}}/\sigma_{\mathbb{N}}^{\alpha}$, with

$$\beta_{\mathbb{N}} = \frac{1}{n^d} \int_{\mathbb{R}^d} \left(\theta \int_{K_n} \psi(x+t)\mathbb{N}(dt) + \sqrt{2\gamma} \int_{K_n} \psi(x+t)G(t)\mathbb{N}(dt)\right)^{\langle\alpha\rangle} dx,$$

and location $1_{\{\alpha=1\}}\mu_{\mathbb{N}}$, with

$$\mu_{\mathbb{N}} = \frac{-2\beta}{\pi} \frac{1}{n^d} \int_{\mathbb{R}^d} \left(\theta \int_{K_n} \psi(x+t)\mathbb{N}(dt) + \sqrt{2\gamma} \int_{K_n} \psi(x+t)G(t)\mathbb{N}(dt)\right)$$
$$\times \log\left(\theta \int_{K_n} \psi(x+t)\mathbb{N}(dt)\right.$$
$$\left. + \sqrt{2\gamma} \int_{K_n} \psi(x+t)G(t)\mathbb{N}(dt)\right) dx.$$

Hence, to determine convergence, we will establish the limits in probability of each of these parameters and thereby determine the joint Fourier–Laplace transform of the sample mean and sample variance. We provide an explicit proof of the convergence of $\sigma_{\mathbb{N}}^{\alpha}$; the proofs for the other two parameters are similar. Let us write

$$H_x = \theta \int_K \psi(x+s)\mathbb{N}(ds) + \sqrt{2\gamma} \int_K \psi(x+s)G(s)\mathbb{N}(ds)$$

so that $\sigma_{\mathbb{N}}^{\alpha} = n^{-d} \int_{\mathbb{R}^d} |H_x|^{\alpha} dx$. Now, $\{H_x\}$ is, conditional on $\mathbb{N}$, a Gaussian process with mean $\theta \int_K \psi(x+s)\mathbb{N}(ds)$ and covariance

$$\text{Cov}_{\mathbb{N}}[H_x, H_y] = \mathbb{E}[(H_x - \mathbb{E}H_x)(H_y - \mathbb{E}H_y)|\mathbb{N}] = 2\gamma \int_K \psi(x+s)\psi(y+s)\mathbb{N}(ds).$$

Taking a second expectation shows that

$$\text{Cov}[H_x, H_y] = \mathbb{E}[(H_x - \mathbb{E}H_x)(H_y - \mathbb{E}H_y)] = 2\gamma \int_K \psi(x+s)\psi(y+s)\,ds,$$

which is zero if $|x_i - y_i| \geq l$ for at least one $i$ between 1 and $d$. Hence, these variables are $l$-dependent (taking the $\infty$-norm for $\mathbb{R}^d$), which will be useful in establishing a weak law of large numbers. It is also true that all moments of $H_x$ exist (even as $n$ increases):

$$|H_x| \leq |\theta| \int_{K_n} |\psi(x+s)|\mathbb{N}(ds) + \sqrt{2\gamma} \int_{K_n} |\psi(x+s)||G(s)|\mathbb{N}(ds)$$
$$\overset{\mathcal{L}}{=} |\theta| \int_{K_n} |\psi(x+s)|\mathbb{N}(ds) + \sqrt{2\gamma}\sqrt{\int_{K_n} \psi(x+s)^2\mathbb{N}(ds)}|G_x|$$
$$\leq \int_{\mathbb{R}^d} |\psi(x+s)|\mathbb{N}(ds)\left(|\theta| + \sqrt{2\gamma}|G_x|\right),$$

where $G_x$ is a dependent sequence of standard normal random variables. The equality in distribution follows from the stability property of Gaussian random variables and the fact that integration with respect to $\mathbb{N}$ is, conditional on $\mathbb{N}$, a discrete sum. The final random variable is Poisson with mean $\int_{\mathbb{R}^d} |\psi(s)| \mathbb{N}(ds)$ multiplied by an independent Gaussian with mean $|\theta|$ and variance $2\gamma$; all moments therefore exist, even as $n \to \infty$. Next, we write

$$\sigma_{\mathbb{N}}^\alpha = \frac{1}{n^d} \int_{\mathbb{R}^d} |H_x|^\alpha \, dx = \frac{1}{n^d} \sum_{j \in \mathbb{Z}^d} \int_B |H_{x+j}|^\alpha \, dx$$

using the same notation as in the proof of Theorem 1. By the same arguments, up to terms going to zero in probability, this expression is the same as $\frac{1}{n^d} \sum_{k \in K_n} \int_B |H_{x-k}|^\alpha \, dx$. Because the variables $J_{-k} = \int_B |H_{x-k}|^\alpha \, dx$ are a random field in $k$ with finite dependence, the weak law of large numbers applies. Hence, $\sigma_{\mathbb{N}}^\alpha = o_P(1) + \frac{1}{n^d} \sum_{k \in K_n} \mathbb{E}[J_{-k}]$ as $n \to \infty$. It remains to compute the expectations. Let

$$\tilde{H}_x = \theta \int_{\mathbb{R}^d} \psi(x+s) \mathbb{N}(ds) + \sqrt{2\gamma} \int_{\mathbb{R}^d} \psi(x+s) G(s) \mathbb{N}(ds)$$

so that

$$\mathbb{E}|\tilde{H}_{x-k}|^\alpha = \mathbb{E}\left| \theta \int_{\mathbb{R}^d} \psi(y) \mathbb{N}(dy) + \sqrt{2\gamma} \sqrt{\int_{\mathbb{R}^d} \psi^2(y) \mathbb{N}(dy)} \tilde{G}_{x-k} \right|^\alpha,$$

where $\tilde{G}_x$ is a mean-zero Gaussian sequence with known correlation structure

$$\mathbb{E}[\tilde{G}_x \tilde{G}_{x+h}] = \frac{\int_{\mathbb{R}^d} \psi(y) \psi(y+h) \, dy}{\int_{\mathbb{R}^d} \psi^2(y) \, dy}.$$

We claim that the average $\overline{\mathbb{E}[J]_n}$ is asymptotically $\mathbb{E}|\tilde{H}_0|^\alpha$. The absolute difference is

$$(19) \quad \begin{aligned} &\left| \frac{1}{n^d} \sum_{k \in K_n} \mathbb{E}[J_{-k}] - \frac{1}{n^d} \sum_{k \in K_n} \mathbb{E} \int_B |\tilde{H}_{x-k}|^\alpha \, dx \right| \\ &\leq \frac{1}{n^d} \sum_{k \in K_n} \int_B \mathbb{E}||H_{x-k}|^\alpha - |\tilde{H}_{x-k}|^\alpha| \, dx. \end{aligned}$$

We have the following inequality, stated as a separate lemma.

LEMMA 1.   *If* $0 < \alpha \leq 1$, *then*

$$||a|^\alpha - |b|^\alpha| \leq |a-b|^\alpha$$

*and if* $1 < \alpha \leq 2$, *then*

$$||a|^\alpha - |b|^\alpha| \leq |a-b|^\alpha + 2\max\{|a|, |b|\}|a-b|^{\alpha/2}$$

*for all real numbers* $a$ *and* $b$.

Proof.    The case $\alpha \leq 1$ is well known. So, suppose that $\alpha > 1$. If $|a| > |b|$, then

$$
\begin{aligned}
|a - b|^\alpha &= (|a - b|^{\alpha/2})^2 \\
&\geq (|a|^{\alpha/2} - |b|^{\alpha/2})^2 \\
&= (|a|^{\alpha/2} - |b|^{\alpha/2})(|a|^{\alpha/2} + |b|^{\alpha/2} - 2|b|^{\alpha/2}) \\
&= (|a|^\alpha - |b|^\alpha) - 2|b|^{\alpha/2}(|a|^{\alpha/2} - |b|^{\alpha/2}),
\end{aligned}
$$

where the second line follows from $\alpha \leq 2$. This, in turn, implies that

$$
\begin{aligned}
|a|^\alpha - |b|^\alpha &\leq |a - b|^\alpha + 2|b|^{\alpha/2}(|a|^{\alpha/2} - |b|^{\alpha/2}) \\
&\leq |a - b|^\alpha + 2\max\{|a|^{\alpha/2}, |b|^{\alpha/2}\}|a - b|^{\alpha/2}.
\end{aligned}
$$

Now, the case $|b| > |a|$ is similar, which proves the lemma.    $\square$

Using this lemma, it suffices to examine the average expected integral of

$$
\begin{aligned}
|H_{x-k} - \tilde{H}_{x-k}|^\delta \\
\leq 2^\delta \left( \left| \theta \int_{K_n^c} \psi(x - k + s)\mathbb{N}(ds) \right|^\delta + \left| \sqrt{2\gamma} \int_{K_n^c} \psi(x - k + s)G(s)\mathbb{N}(ds) \right|^\delta \right),
\end{aligned}
$$

where $\delta$ is either $\alpha$ or $\alpha/2$. If $\delta = \alpha/2$, we have

$$
\begin{aligned}
\frac{1}{n^d} \sum_{k \in K_n} \int_B \mathbb{E}|H_{x-k} - \tilde{H}_{x-k}|^{\alpha/2} \\
\leq 2^{\alpha/2} \frac{1}{n^d} \sum_{k \in K_n} \int_B \int_{K_n^c} |\psi(x - k + s)|^{\alpha/2} \, ds \, dx \, (|\theta|^{\alpha/2} + (2\gamma)^{\alpha/4}),
\end{aligned}
$$

using the fact that $\mathbb{E}\mathbb{N}(ds) = ds$. If $\delta = \alpha$, using the result that

$$
\begin{aligned}
\mathbb{E} \left| \theta \int_{K_n^c} \psi(x - k + s)\mathbb{N}(ds) \right|^\alpha + \mathbb{E} \left| \sqrt{2\gamma} \int_{K_n^c} \psi(x - k + s)G(s)\mathbb{N}(ds) \right|^\alpha \\
\leq |\theta|^\alpha \mathbb{E} \left( \int_{K_n^c} |\psi(x - k + s)|^{\alpha/2}\mathbb{N}(ds) \right)^2 \\
+ (2\gamma)^{\alpha/2} \mathbb{E} \left( \int_{K_n^c} |\psi(x - k + s)|^{\alpha/2}|G(s)|^{\alpha/2}\mathbb{N}(ds) \right)^2 \\
= |\theta|^\alpha \int_{K_n^c} |\psi(x - k + s)|^\alpha \, ds \\
+ (2\gamma)^{\alpha/2} \Bigg[ \int_{K_n^c} |\psi(x - k + s)|^\alpha \, ds \\
+ \left( \int_{K_n^c} |\psi(x - k + s)|^{\alpha/2} \, ds \right)^2 \Bigg] \mathbb{E}|G|^\alpha,
\end{aligned}
$$

we easily obtain

$$\frac{1}{n^d} \sum_{k \in K_n} \int_B \mathbb{E}|H_{x-k} - \tilde{H}_{x-k}|^\alpha$$

$$\leq 2^\alpha \frac{1}{n^d} \sum_{k \in K_n} \int_B \int_{K_n^c} |\psi(x-k+s)|^\alpha \, ds \, dx \, (|\theta|^\alpha + (2\gamma)^{\alpha/2} \mathbb{E}|G|^\alpha)$$

$$+ 2^\alpha \frac{1}{n^d} \sum_{k \in K_n} \int_B \left( \int_{K_n^c} |\psi(x-k+s)|^{\alpha/2} \, ds \right)^2 dx \, ((2\gamma)^{\alpha/2} \mathbb{E}|G|^\alpha).$$

So, in order to show that (19) tends to zero as $n \to \infty$, it is enough to demonstrate that

$$\frac{1}{n^d} \sum_{k \in K_n} \int_B \left( \int_{K_n^c} |\psi(x-k+s)|^\delta \, ds \right)^\phi dx$$

tends to zero, for $\delta$ equal to either $\alpha$ or $\alpha/2$ and $\phi$ equal to 1 or 2. Note that through a simple change of variable, this becomes

$$\frac{1}{n^d} \int_{K_n} \left( \int_{K_n^c} |\psi(s-x)|^\delta \, ds \right)^\phi dx.$$

It is a simple but tedious analysis exercise to show that this tends to zero as $n \to \infty$, for $\delta = \alpha$ or $\alpha/2$ and $\phi$ equal to 1 or 2. The result of this analysis is that

$$\sigma_\mathbb{N}^\alpha \xrightarrow{P} \lambda(K) \mathbb{E} \left| \theta \int_{\mathbb{R}^d} \psi(y) \mathbb{N}(dy) + \sqrt{2\gamma} \sqrt{\int_{\mathbb{R}^d} \psi^2(y) \mathbb{N}(dy)} Z \right|^\alpha = \lambda(K) \tilde{\sigma}_\infty^\alpha(\theta, \gamma),$$

where $Z$ is a standard normal random variable. Note that this limiting scale parameter is not random. Using similar techniques, one can show that $\beta_\mathbb{N}$ converges in probability to

$$\lambda(K) \mathbb{E} \left( \theta \int_{\mathbb{R}^d} \psi(y) \mathbb{N}(dy) + \sqrt{2\gamma} \sqrt{\int_{\mathbb{R}^d} \psi^2(y) \mathbb{N}(dy)} Z \right)^{\langle \alpha \rangle} = \lambda(K) \tilde{\beta}_\infty(\theta, \gamma)$$

and that $\mu_\mathbb{N}$ tends to

$$\lambda(K) \frac{-2\beta}{\pi} \mathbb{E} \Bigg[ \left( \theta \int_{\mathbb{R}^d} \psi(y) \mathbb{N}(dy) + \sqrt{2\gamma} \sqrt{\int_{\mathbb{R}^d} \psi^2(y) \mathbb{N}(dy)} Z \right)$$

$$\times \log \left( \theta \int_{\mathbb{R}^d} \psi(y) \mathbb{N}(dy) + \sqrt{2\gamma} \sqrt{\int_{\mathbb{R}^d} \psi^2(y) \mathbb{N}(dy)} Z \right) \Bigg],$$

which is $\lambda(K)\tilde{\mu}_\infty(\theta,\gamma)$. Now, multiplying our results by $r^{-1/\alpha}\lambda(K)^{-1/\alpha}$, our joint Fourier–Laplace transform is

$$
(20)
\begin{aligned}
&\mathbb{E}\exp\left\{-r^{-1}\lambda(K)^{-1}\sigma_\mathbb{N}^\alpha\left(1-i\beta\frac{\beta_\mathbb{N}}{\sigma_\mathbb{N}^\alpha}\right)+ir^{-1}\lambda(K)^{-1}\mu_\mathbb{N}1_{\{\alpha=1\}}\right\}\\
&\qquad\to\mathbb{E}\exp\left\{-r^{-1}\tilde{\sigma}_\infty^\alpha(\theta,\gamma)\left(1-i\frac{\tilde{\beta}_\infty(\theta,\gamma)}{\tilde{\sigma}_\infty^\alpha(\theta,\gamma)}\right)+ir^{-1}\tilde{\mu}_\infty(\theta,\gamma)1_{\{\alpha=1\}}\right\},
\end{aligned}
$$

by the dominated convergence theorem. The existence of $\tilde{S}_\infty(\alpha)$ and $\tilde{U}_\infty(\alpha)$ now follows from the continuity of the limiting Fourier–Laplace transform at $(0,0)$. Letting $\gamma=0$, one recognizes the Fourier transform of an $\alpha$-stable random variable with parameters as described in Theorem 2, which is $\tilde{S}_\infty(\alpha)$; if $\theta=0$, then one obtains the Laplace transform of a positive $\alpha/2$-stable random variable $\tilde{U}_\infty(\alpha)$, whose scale parameter is calculated in the statement of the theorem.

Finally, we must remove the truncation of the model. Define

$$
X_l(t)=\int_{\mathbb{R}^d}\psi_l(x+t)\mathbb{M}(dx)=\int_{\mathbb{R}^d}\psi(x+t)1_{lD}(x+t)\mathbb{M}(dx),
$$

where $D=(-1,1]^d$ so that $lD$ is a $d$-dimensional cube centered at the origin with width $2l$. Then

$$
\begin{aligned}
&n^{-d/\alpha}\int_{K_n}X(t)\mathbb{N}(dt)-n^{-d/\alpha}\int_{K_n}X_m(t)\mathbb{N}(dt)\\
&\qquad=n^{-d/\alpha}\int_{\mathbb{R}^d}\int_{K_n}\psi(x+t)1_{\{lD\}^c}(x+t)\mathbb{N}(dt)\mathbb{M}(dx)
\end{aligned}
$$

must tend to zero in probability as $l\to\infty$ for any fixed $n$. Taking the $\alpha$ power of the scale of the above random variable, conditional on $\mathbb{N}$, we obtain

$$
\frac{1}{n^d}\int_{\mathbb{R}^d}\left|\int_{K_n}\psi(x+t)1_{\{lD\}^c}(x+t)\mathbb{N}(dt)\right|^\alpha\lambda(dx).
$$

If $\alpha<1$, this is bounded by

$$
\begin{aligned}
&\frac{1}{n^d}\int_{\mathbb{R}^d}\int_{K_n}|\psi(x+t)|^\alpha 1_{\{lD\}^c}(x+t)\mathbb{N}(dt)\lambda(dx)\\
&\qquad=\frac{1}{n^d}\int_{K_n}\int_{\{lD\}^c-t}|\psi(x+t)|^\alpha\lambda(dx)\mathbb{N}(dt)\\
&\qquad=\int_{\{lD\}^c}|\psi(y)|^\alpha\lambda(dy)\frac{\mathbb{N}(K_n)}{n^d}.
\end{aligned}
$$

Thus, the limit superior as $n\to\infty$ of $\frac{1}{n^d}\int_{\mathbb{R}^d}|\int_{K_n}\psi(x+t)1_{\{lD\}^c}(x+t)\mathbb{N}(dt)|^\alpha\times\lambda(dx)$ tends to zero as $l\to\infty$. When $\alpha>1$, we use the bound of

$$\sup_x \left| \int_{K_n} \psi(x+t) 1_{\{lD\}^c}(x+t) \mathbb{N}(dt) \right|$$

$$\times \frac{1}{n^d} \int_{\mathbb{R}^d} \left| \int_{K_n} \psi(x+t) 1_{\{lD\}^c}(x+t) \mathbb{N}(dt) \right|^{\alpha-1} \lambda(dx)$$

$$\leq \int_{\mathbb{R}^d} |\psi(t)| \mathbb{N}(dt) \cdot \frac{1}{n^d} \int_{\mathbb{R}^d} \int_{K_n} |\psi(x+t)|^{\alpha-1} 1_{\{lD\}^c}(x+t) \mathbb{N}(dt) \lambda(dx)$$

$$= Q \cdot \int_{\{lD\}^c} |\psi(y)|^{\alpha-1} \, dy,$$

where $Q = \int_{\mathbb{R}^d} |\psi(t)| \mathbb{N}(dt)$ is a positive random variable that is bounded in $n$. So the limit superior in this case also tends to zero as $l \to \infty$. Similar arguments can be applied to the sum of squares and since the limiting joint Fourier–Laplace transform is continuous in $l$, we can take the limit as $l \to \infty$ on both sides of our joint weak convergence (20). This completes the proof. $\square$

PROOF OF THEOREM 3.  This proof has the same structure as that of Theorem 6.3.1 from PRW [8]. First, we show that $\tau_{\mathbb{N}(B+y)}$ is almost surely asymptotic to $\tau_{\Lambda(B+y)} = \tau_{\Lambda(B)}$. As in the proof of Theorem 2, $\frac{\mathbb{N}(B+y)}{\Lambda(B+y)} \xrightarrow{\text{a.s.}} 1$ as $n \to \infty$; this uses the condition that $c_n \delta(K_n) \to \infty$. It follows from the form of $\tau(u)$ that $\frac{\tau_{\mathbb{N}(B+y)}}{\tau_{\Lambda(B+y)}} \xrightarrow{\text{a.s.}} 1$. Let $x$ be a continuity point of $J(x)$, the cumulative distribution function of the limit random variable $J$. Then

$$\left\{ \tau_{\mathbb{N}(B+y)} \left( \frac{\hat{\mu}_{K_n, B, y} - \hat{\mu}_{K_n}}{\hat{\sigma}_{K_n, B, y}} \right) \leq x \right\}$$

$$= \left\{ \tau_{\mathbb{N}(B+y)} \left( \frac{\hat{\mu}_{K_n, B, y} - \mu}{\hat{\sigma}_{K_n, B, y}} \right) \leq x + \tau_{\mathbb{N}(B+y)} \left( \frac{\hat{\mu}_{K_n} - \mu}{\hat{\sigma}_{K_n, B, y}} \right) \right\}.$$

For any $t > 0$, let

$$R_{K_n, B}(t) = \lambda(K_n(1-c))^{-1} \int_{K_n(1-c)} 1_{\{\tau_{\mathbb{N}(B+y)}(\frac{\hat{\mu}_{K_n} - \mu}{\hat{\sigma}_{K_n, B, y}}) \leq t\}} \, dy$$

$$= \lambda(K_n(1-c))^{-1} \int_{K_n(1-c)} 1_{\{d_{\mathbb{N}(B+y)} \hat{\sigma}_{K_n, B, y} \geq a_{\mathbb{N}(B+y)}(\hat{\mu}_{K_n} - \mu)/t\}} \, dy.$$

Now, for all $\delta > 0$, $a_{\mathbb{N}(B+y)}(\hat{\mu}_{K_n} - \mu) \leq \delta$ with probability tending to one, since $a_{\mathbb{N}(B+y)}/a_{\mathbb{N}(K_n)} \xrightarrow{P} 0$ follows from $c_n \delta(K_n) \to \infty$. So, with probability tending to one,

$$R_{K_n, B}(t) \geq \lambda(K_n(1-c))^{-1} \int_{K_n(1-c)} 1_{\{d_{\mathbb{N}(B+y)} \hat{\sigma}_{K_n, B, y} \geq \delta/t\}} \, dy.$$

If $\delta/t$ is a continuity point of $W$, then the above expression tends to $\mathbb{P}[W \geq \delta/t]$, by Theorem 6.3.1 of PRW [8]. Since $W$ has no point mass at zero, we can make $R_{K_n,B}(t)$ arbitrarily close to 1 by choosing $\delta$ sufficiently small. So, for all $t$, $R_{K_n,B}(t) \xrightarrow{P} 1$. Next, using the inequality

$$1_{\{\tau_{\mathbb{N}(B+y)}(\frac{\hat{\mu}_{K_n,B,y}-\mu}{\hat{\sigma}_{K_n,B,y}}) \leq x + \tau_{\mathbb{N}(B+y)}(\frac{\hat{\mu}_{K_n}-\mu}{\hat{\sigma}_{K_n,B,y}})\}}$$
$$\leq 1_{\{\tau_{\mathbb{N}(B+y)}(\frac{\hat{\mu}_{K_n,B,y}-\mu}{\hat{\sigma}_{K_n,B,y}}) \leq x+t\}} + 1_{\{\tau_{\mathbb{N}(B+y)}(\frac{\hat{\mu}_{K_n}-\mu}{\hat{\sigma}_{K_n,B,y}}) > t\}},$$

we can establish that

$$L_{K_n,B}(x)$$
$$\leq \lambda(K_n(1-c))^{-1} \int_{K_n(1-c)} 1_{\{\tau_{\mathbb{N}(B+y)}(\frac{\hat{\mu}_{K_n,B,y}-\mu}{\hat{\sigma}_{K_n,B,y}}) \leq x+t\}} \, dy + (1 - R_{K_n,B}(t)).$$

Now, by the first convergence in (16), we may apply Theorem 6.3.1 of PRW [8] to

$$\tau_{\mathbb{N}(B+y)}\left(\frac{\hat{\mu}_{K_n,B,y}-\mu}{\hat{\sigma}_{K_n,B,y}}\right)$$

and obtain, for any $\epsilon > 0$, $L_{K_n,B}(x) \leq J(x+t) + \epsilon$ with probability tending to one. At this point, let $t$ tend to zero. Similar arguments produce the opposite inequality $L_{K_n,B}(x) \geq J(x+t) - \epsilon$. Now letting $\epsilon \to 0$, we obtain $L_{K_n,B}(x) \xrightarrow{P} J(x)$, as desired. The proofs of (ii) and (iii) are similar to the proof of Theorem 2.2.1 in PRW [8]. $\quad\square$

U.S. BUREAU OF THE CENSUS
STATISTICAL RESEARCH DIVISION
4700 SILVER HILL ROAD
WASHINGTON, DISTRICT OF COLUMBIA 20233-9100
USA
E-MAIL: tucker.s.mcelroy@census.gov

DEPARTMENT OF MATHEMATICS
UNIVERSITY OF CALIFORNIA, SAN DIEGO
LA JOLLA, CALIFORNIA 92093-0112
USA
E-MAIL: politis@math.ucsd.edu